\documentclass[ a4paper, draft,10pt]{article}

\usepackage{amsmath, amssymb, amsfonts}

\usepackage[T1]{fontenc}
\usepackage[cp1250]{inputenc}

\newtheorem{thm}{Theorem}[section]
\newtheorem{lem}[thm]{Lemma}

\newtheorem{stm}[thm]{Statement}

\newtheorem{rmk}{Remark}

\title{On estimation of solutions of neutral type systems on nonclosed sets\footnote{This work was partially supported by Polish National Science Centre grant No. N N514 238438.}}
\author{G.M. Sklyar, P. Polak
\\
sklar@univ.szczecin.pl, piotr.polak@wmf.univ.szczecin.pl}

\date{Institute of Mathematics, University of Szczecin, Wielkopolska 15, 70-451 Szczecin, Poland}

\begin{document}
\maketitle




\abstract{We consider a differential system of neutral type with distributed delay.
We obtain a precise norm estimation of solutions of the system in question on some nonclosed set. Our result is based on a spectral analysis of the operator and Riesz basis theory. We also discuss the stability properties of corresponding solutions.


\section{Introduction}
The important problem in the theory of differential equations is to determine the asymptotic behavior of the solutions. Even in the case 
when solutions tend to the equilibrium it may be important to see how fast they approach it. If the equilibrium is zero then those questions 
relate directly to the stability analysis of a system. In \cite{SkPo13} we gave an estimation of growth speed of all solutions of neutral type systems of some class. Recently the estimations of the solutions growth which are not valid on whole state space but for initial states from some 
nonclosed set were described. A. Borichev and Y. Tomilov characterized in \cite{BoTo10} the decay rate of solutions of an abstract linear differential
equation in the Hilbert space when initial states were from a certain nonclosed dense set. Their characterization was given in resolvent 
terms. This encouraged us to study the problem of the rate of growth (or decay) of solutions of neutral type equations in the same context. We 
give an extension of those result in the case of neutral type equations without an assumption about boundedness of corresponding semigroup. Our 
concept is based on the theorems on maximal asymptotics \cite{Sk10,Sk10a} and our previous results concerning some class of delay systems of neutral type \cite{SkPo13}.

\section{Preliminaries}
In the present work, following \cite{RSR05}, we consider the delay systems of neutral type of the form
\begin{equation}\label{rowzopu}
 \dot{z}(t)=A_{-1}\dot{z}(t-1) + \int_{-1}^0A_2(\theta)\dot{z}(t+\theta)d\theta 
+ \int_{-1}^0 A_3(\theta)z(t+\theta)d\theta,
\end{equation}
where $A_{-1}$ is a $n \times n$ invertible complex matrix, $A_2$ and $A_3$ are $n \times n$ matrices of functions from $L_2(-1,0)$. 
We rewrite equation (\ref{rowzopu}) in the operator form
\begin{equation}\label{operatorform}
\dot{x}=\mathcal{A}x, \qquad x \in M_2,
\end{equation}
where $M_2=\mathbb{C}^n \times L_2(-1,0;\mathbb{C}^n)$, the operator $\mathcal{A}$ is then given by
\begin{equation}\label{operatorA}
 \mathcal{A} \binom{y(t)}{z_t(\cdot)} =
\binom{\int_{-1}^0A_2(\theta)\dot{z}_t(\theta)d\theta +
 \int_{-1}^0 A_3(\theta)z_t(\theta)d\theta}{dz_t(\theta)/d\theta}, 
\end{equation}
and the domain of $\mathcal{A}$ is as follows:
\begin{equation}\label{domain}
 \mathcal{D}(\mathcal{A})= \left\{ (y,z(\cdot)):z \in H^1(-1,0;\mathbb{C}^n), y=z(0)-A_{-1}z(-1) \right\}\subset M_2.
\end{equation}

 Denote the eigenvalues of the matrix $A_{-1}$ by $\mu_m, m=1,\ldots, \ell$ ($|\mu_1|\geq|\mu_2|\geq \ldots \geq |\mu_{\ell}|$), and their multiplicities by $p_m$ ($ \sum{p_m}=n$). Without loss of generality we assume that if $|\mu_1|=|\mu_2|=\ldots =|\mu_{\ell_0}|$ then $p_1\geq p_2\geq \ldots \geq p_{\ell_0}$.
  The eigenvalues of $\tilde{\mathcal{A}}$ are complex logarithms of $\mu_m$ and zero i.e. 
\begin{eqnarray*} \sigma(\tilde{\mathcal{A}})&\!\!\!\!=& \!\!\!\!\{ \tilde{\lambda}_m^{(k)}=\ln|\mu_m|+i(\arg\mu_m+2k\pi),\mu_m\in \sigma(A_{-1}), m=1,\ldots,\ell  ;k \in \mathbb{Z} \}\\&&\cup \{0\}.\end{eqnarray*}
Almost all eigenvalues of $ \mathcal{A}$ lie close to $\tilde{\lambda}_m^{(k)}$. More precisely, for $k$ large enough they are contained in the discs $L_m^{(k)}$ centered at $\tilde{\lambda}_m^{(k)}$ of radii $r_k \to 0$ (see Theorem 4 \cite{RSR08}).
The sum of multiplicities of eigenvalues of $ \mathcal{A}$ lying in each disc centered at $\tilde{\lambda}_m^{(k)}$ equals the multiplicity of $\tilde{\lambda}_m^{(k)}$ and $\mu_m$, that is $p_m$. 
We denote eigenvalues of operator $\mathcal{A}$ by $\lambda_{m,i}^{(k)}, k\in \mathbb{Z}; m=1,\ldots,\ell$, and we have $\{\lambda_{m,i}^{(k)}\}_{i=1}^{p_m}\subset L_m^{(k)}, |k|>N; m=1,\ldots,\ell$.
Let us denote $\mathcal{A}$-invariant subspaces $V_m^{(k)}=P_m^{(k)}M_2$, where  $P_m^{(k)}x=\frac1{2\pi i} \int_{L_m^{(k)}} R(\mathcal{A},\lambda)xd\lambda$ are Riesz projectors, $m=1,\ldots,\ell, k \in \mathbb{Z}$.
The sequence of $p_m$-dimensional subspaces $V_m^{(k)}, m=1,\ldots,\ell, |k| \geq N $, and some $2(N+1)n$-dimensional subspace $W_N$ constitute $\mathcal{A}$-invariant Riesz basis of space $M_2$. 
Notice that  
the subspaces $\tilde{V}_m^{(k)}, m=1,\ldots,\ell, |k| \geq N $, and some $2(N+1)n$-dimensional subspace $\tilde{W}_N$ constitute $\tilde{\mathcal{A}}$-invariant Riesz basis of the space $M_2$, which is quadratically close to the previous one. We denote $ \mathcal{A}|_{V_m^{(k)}}$ by $\mathcal{A}_m^{(k)}$, for each $m=1,\ldots, \ell; |k|>N$ and $\mathcal{A}|_{W_N}$ by $\mathcal{A}^{(N)}$. 
The same for operator $\tilde{\mathcal{A}}$.  
Theorem 7 states in particular that to each $\tilde{\lambda}_m^{(k)}\in \sigma(\tilde{\mathcal{A}}) \setminus \{0\}$ and each Jordan chain of eigen- and rootvectors of the matrix $A_{-1}$ corresponds the Jordan chain of $\tilde{\mathcal{A}}:v_{m,j}^{(k),0}, v_{m,j}^{(k),1}, \ldots, v_{m,j}^{(k),p_{m,j}-1},$ i.e. the vectors $v^s$ (the indices $k,m,j$ are omitted) satisfy the relation $(\tilde{\mathcal{A}}-\lambda I)v^s=v^{s-1}$. 
They are given by 
\begin{equation}\label{eigenvectors}
v_{m,j}^{(k),s}=
\binom{0}{e^{\tilde{\lambda}_m^{(k)} \theta}P^s_{m,j}(\theta)},
\end{equation}
where $P^s_{m,j}(\theta)$ is some polynomial independent of $k$, $s=0,\ldots, p_{m,j}-1; m=1,\ldots,\ell; k\in \mathbb{Z};j=1,\ldots,\nu_m$.
Some of the consequences of this characterization for the norm of semigroup $e^{\tilde{\mathcal{A}}t}$ are given in \cite{SkPo13}

\section{Growth bound estimation}  
With assumption that for all eigenvalues of operator $\mathcal{A}$ their real parts lie at the left side of $\tilde{\omega}$ it was proven that 
\begin{equation}\label{OldEstim}
\|e^{\mathcal{A}t}\| \leq M_p e^{\omega t}(t^{p-1}), \quad t\geq 1, 
\end{equation}
where $p$ is the sum of sizes of all Jordan blocks corresponding to the maximal eigenvalue of $A_{-1}$. Such assumption implies that there exist a sequence of eigenvalues of $\mathcal{A}$ whose real parts tend to $\tilde{\omega}$. If this tending is arbitrary slow or initial state is sufficiently smooth then the power of $t$ in the estimation \eqref{OldEstim} can be decreased. Namely we have the following
%
\begin{thm}\label{thm1}
We consider system \eqref{operatorform}. If $\rm{Re}\,\lambda < \omega$ for all $\lambda \in \sigma(\mathcal{A})$ and \, $\rm{Re}\lambda_{1,i}^{(k)} \leq \omega-|k|^{-s}, \, k\in \mathbb{Z}, i=1,\ldots,p_1 $, where $\{\lambda_{1,i}^{(k)}\}_{i=1}^{p_1} \subset L_{1}^{(k)}$ then 
$$\|e^{\mathcal{A}t}\mathcal{A}^{-n}x\| \leq Ce^{\omega t}t^{p-1-\frac ns}\|x\|,$$
\end{thm}
Before we prove Theorem \ref{thm1} we give two lemmas
\begin{lem}\label{lem1}
We consider system \eqref{operatorform}. If $\rm{Re}\,\lambda < \omega$ for all $\lambda \in \sigma(\mathcal{A})$ and for some $s\in \mathbb{N}$ holds $\rm{Re}\lambda_{1,i}^{(k)} \leq \omega-|k|^{-s}, \, k\in \mathbb{Z}, i=1,\ldots,p_1 $, where $\{\lambda_{1,i}^{(k)}\}_{i=1}^{p_1} \subset L_{1}^{(k)}$ then there exist a constant $C_1>0$ independent of $k$ such that
$$\|e^{\mathcal{A}t}|_{V_1^{(k)}}\| \leq C_1e^{(\omega-|k|^{-s}) t}t^{p-1}, \quad t>1.$$
\end{lem}
{\it Proof} Because the families of subspaces $\{V_m^{(k)}\}_{k\in \mathbb{Z};m=1,\ldots,\ell}$ and $\{\tilde{V}_m^{(k)}\}_{k\in \mathbb{Z};m=1,\ldots,\ell}$ constitute quadratically close Riesz basis, there exists a bounded operator $T_N$, with bounded inverse, which is close to identity and transforms almost all subspaces $V_m^{(k)}$ to $\tilde{V}_m^{(k)}$. 
Such operator $T_N$ can be defined on every basis subspace $V_m^{(k)}$ by the formula
\begin{equation}\label{tn1}
T_N|_{V_m^{(k)}}x=\tilde{P}_m^{(k)}x, \,\mbox{for } |k|>N, \, m=1,2,\ldots,\ell,
\end{equation}
and $T_N|_{W_N}x=x$.
It is easy to see that the operator $T_N$ is bounded on $M_2$ and close to identity. Therefore $T_N$ is invertible, its inverse $T_N^{-1}$ is bounded and transforms all but finitely many  subspaces $V_m^{(k)}$ onto $\tilde{V}_m^{(k)}$. 
Now we denote $\mathcal{A}_k^0=\mathcal{A}_1^{(k)}-\tilde{\lambda}_1^{(k)}I, |k|>N$
and using $T_N$ we define the operator $\mathcal{B}^0_k:\tilde{V}_k \to \tilde{V}_k$ close to $\mathcal{A}_k^0:V_k \to V_k$ for each $|k|>N$ by the formula $\mathcal{B}^0_k=T_N\mathcal{A}_k^0T_N^{-1}$. Let $\tilde{x}=T_Nx$.
Then we get 
\begin{equation}
 \|e^{ \mathcal{A}_1^{(k)}t}\|=|e^{\tilde{\lambda}_1^{(k)}t}|\cdot \|T_N^{-1}e^{\mathcal{B}^0_kt}T_N\|, 
\end{equation}
which can be rewritten as  
\begin{equation}\label{BkA1}
 \|e^{\mathcal{A}_1^{(k)}t}\|\leq Ce^{\omega t}\|e^{\mathcal{B}^0_kt}\|, 
\end{equation}
where $C$ is a positive constant.
It was shown in Lemma 2.2 \cite{SkPo13} that for each $\delta>0$ there exist $k_0$ large enough such that 
\begin{equation}\label{BkA0}
\|\mathcal{B}^0_k-\tilde{\mathcal{A}}_k^0\| \leq \delta \mbox{ for any } |k|>k_0,
\end{equation} 
where $\tilde{\mathcal{A}}_k^0=\tilde{\mathcal{A}}_1^{(k)}-\tilde{\lambda}_1^{(k)}$.

The eigen- and rootvectors of operator $\tilde{\mathcal{A}}$ are given by \eqref{eigenvectors}, we can see that 

the number and lengths of all Jordan chains of operators $\tilde{\mathcal{A}}_k=\mathcal{\tilde A} |_{\tilde V_k}$ are independent of $k$. Thus 
all operators $\tilde{\mathcal{A}}_k^0:=\tilde{\mathcal{A}}_k-\tilde{\lambda}_kI$ have the same matrix, say $A_0$ in the basis $\{v^{(k),s}\}_{s=1}^{p_1}$. Moreover if $\tilde{\mathcal{A}}_k^0=S_k^{-1}A_0S_k$, where $S_k:\tilde{V}_1^{(k)} \to \mathbb{C}^{p_1}$ then $\|S_k\|,\|S_k^{-1}\|$ are uniformly bounded. Therefore denoting the matrix of operator $\mathcal{B}^0_k$ by $B^0_k$ we have 
$$
\|e^{\mathcal{B}^0_kt}\| \leq C_1 \|e^{B^0_kt}\|,
$$
which can be rewritten as   
\begin{equation}\label{Bk1}
 \|e^{\mathcal{B}^0_kt}\| \leq C_1 e^{-|k|^st}\|e^{(B^0_k+|k|^sI)t}\|.
\end{equation}
Inequality \eqref{BkA0} is also satisfied for matrices of operators $\mathcal{B}^0_k, \tilde{\mathcal{A}}_k^0$ and it is easy to see that for any $\delta>0$ there exists $k_0$ such that 
$$\|(B^0_k+|k|^sI)-A_0\| \leq \delta,  |k|>k_0.$$ 
Hence we can apply Lemma 2.3 \cite{SkPo13} to the family of matrices $\{B^0_k+|k|^sI\}_{|k|>k_0}$, we get
\begin{equation}\label{Bk2}
\|e^{(B^0_k+|k|^sI)t}\| \leq C_2t^{p_1-1}, t>1, |k|>k_0,
\end{equation}
where we used the fact that $\sigma(B^0_k+|k|^sI) \subset \mathbb{C}^-$.
Combining \eqref{BkA1}, \eqref{Bk1}, \eqref{Bk2} we get 
$$\|e^{\mathcal{A}_1^{(k)}t}\| \leq C_1e^{(\omega-|k|^{-m}) t}t^{p_1-1}, \quad t>1, k>k_0,$$
which ends the proof of Lemma \ref{lem1}.
$\hfil \square\ $

\begin{lem}\label{lem2}
We consider system \eqref{operatorform}. There exist $k_0 \in \mathbb{N}$ and constant $C_2>0$ such that 
$$ \|\mathcal{A}_k^{-1}\|\leq \frac{C_2}{|k|}, \quad |k| \geq k_0. $$
\end{lem}
{\it Proof}
Using  uniformly bounded operators $T_N, S_k$ described in the proof of Lemma \ref{lem1} we define operator $\mathcal{B}_k:\tilde{V}_k \to \tilde{V}_k$ by the formula $\mathcal{B}_k=T_N\mathcal{A}_kT_N^{-1}, |k|>N$. Notice that the matrix of operator $\mathcal{B}_k$ is given by $B_k=S_k^{-1}\mathcal{B}_kS_k$. Hence
\begin{equation*}\label{bk}
\|\mathcal{A}_k^{-1}\|\leq C_3 \|B_k^{-1}\|,
\end{equation*}
for $|k|$ large enough.
Moreover for corresponding matrices of operators $\mathcal{B}_k, \tilde{\mathcal{A}}_k$ and some constant $C_4$ we have $\|B_k-\tilde{A}_k\|\leq C_4 \|\mathcal{B}_k-\tilde{\mathcal{A}}_k\|, |k|>N$.
Lemma 2.2 \cite{SkPo13} implies that there exists $k_0\in \mathbb{N}$ such that $\|\mathcal{B}_k-\tilde{A}_k\|\leq C_4^{-1}$ for $|k|>k_0$, thus we have $\|B_k-\tilde{A}_k\|\leq 1$. 
Matrix $\tilde{A}_k$ consist of Jordan blocks with eigenvalue $\tilde{\lambda}_1^{(k)}=\tilde{\omega}+i(\arg \mu_1+2k\pi)$. 
Therefore we are able to use Statement \ref{stm2} and we obtain
\begin{equation*}\label{bk1}
\|\mathcal{A}_k^{-1}\|\leq \frac{CC_3}{|\tilde{\omega}+i(\arg \mu_1+2k\pi)|} \leq \frac{C_2}{|k|}, \qquad |k|>k_0,
\end{equation*}
for some constant $C_2$. 
$\hfil \square\ $

{\it Proof of Theorem \ref{thm1}} 
Without loss of generality we assume that all eigenvalues of matrix $A_{-1}$ have the same modulus, it means that all eigenvalues of operator $\tilde{\mathcal{A}}$ lie on the line $x=\tilde{\omega}$. Therefore we omit index $m$ and write $\lambda_{i}^{(k)}, L_k, V_k, $ instead of $\lambda_{m,i}^{(k)}, L_m^{(k)}, V_m^{(k)}$.
We decompose $x$ in the Riesz Basis $\{V_k\}_{|k|>N} \cup {W_N}$ i.e. $x=\sum_{|k|\geq N} x_k, x_k\in V_k$ and we have
$$\|e^{\mathcal{A}t}\mathcal{A}^{-n}x\|^2\leq \sum_{|k|\geq N}\|e^{\mathcal{A}_kt}\mathcal{A}_k^{-n}x_k\|^2 \leq \sum_{|k|\geq N}\|e^{\mathcal{A}_kt}\|^2\|\mathcal{A}_k^{-1}\|^{2n}\|x_k\|^{2}.$$
We apply Lemma \ref{lem1} and Lemma \ref{lem2} in above inequality and we obtain 
$$\|e^{\mathcal{A}t}\mathcal{A}^{-n}x\|^2\leq C_1^2C_2^2\sum_{|k|\geq N}e^{2(\omega-|k|^{-s}) t}t^{2(p-1)} |k|^{-2n} \|x_k\|^2.$$
Gathering the terms independent of $k$ and additional term $t^{-2\frac{n}{s}}$ we get
$$\|e^{\mathcal{A}t}\mathcal{A}^{-n}x\|^2\leq \left(C_1 C_2 e^{\omega t} t^{p-1-\frac{n}{s}}\right)^2\sum_{|k|\geq N}\left(e^{-|k|^{-s}t} (|k|^{-s}t)^{\frac{n}{s}}\right)^2 \|x_k\|^2.$$
Since function $f(x)=e^{-x}x^\frac{n}{s}, \,\, x\geq 0$ is bounded by some constant $M>0$, therefore \\
$e^{-|k|^{-s}t} (|k|^{-s}t)^{\frac{n}{s}} \leq M$ for any $k\in \mathbb{Z}$ and we have 
$$\|e^{\mathcal{A}t}\mathcal{A}^{-n}x\|^2\leq \left(M C_1 C_2 e^{\omega t} t^{p-1-\frac{n}{s}}\right)^2\sum_{|k|\geq N} \|x_k\|^2.$$
Subspaces $\{V_k\}_{|k|>N} \cup {W_N}$ constitute Riesz basis, thus there exists constant $C_3>0$ such that $\sum_{|k|\geq N} \|x_k\|^2\leq C_3^2 \|x\|^2$ and we finaly obtain
$$\|e^{\mathcal{A}t}\mathcal{A}^{-n}x\|\leq  C e^{\omega t} t^{p-1-\frac{n}{s}} \|x\|,$$
where $C=M C_1 C_2$ is a new constant. 
$\hfil \square\ $

\section{Appendix}
\begin{stm}\label{stm1}
Let $A_{\lambda}=[a_{i,j}]\in M_n(\mathbb{C})$ be a Jordan block of egenvalue $\lambda$, where $|\lambda|\geq 1$ and $B_{\lambda}=[b_{i,j}]\in M_n(\mathbb{C})\in M_n(\mathbb{C})$ be such that $\|B_{\lambda}-A_{\lambda}\|\leq 1$, where $\|A\|=\sum|a_{i,j}|$, then there exist constant $M$ such that 
$$|\det B_{\lambda}| \leq M |\lambda|^n, $$ 
and for $|\lambda|\geq 2M$ we have also 
$$ |\det B_{\lambda}|\geq \frac12 |\lambda|^n.$$
\end{stm}
{\it Proof} 
Let us define $\varepsilon_{i,j}=b_{i,j}-a_{i,j}$, 
it is easy to see that 
\begin{equation}\label{detbl}
	|\det B_{\lambda}|=|\lambda^n+f_1(\varepsilon_1,\ldots,\varepsilon_{n^2})\lambda^{n-1}+f_2(\varepsilon_1,\ldots,\varepsilon_{n^2})\lambda^{n-2}+\ldots+
	f_n(\varepsilon_1,\ldots,\varepsilon_{n^2})\lambda^{0}|,
\end{equation}
where $f_1,\ldots,f_n$ are polynomials. From the assumption that $\|B_{\lambda}-A_{\lambda}\|\leq 1$ we have  $\sum|\varepsilon_{i,j}|\leq 1$ and therefore 
$$ M_0:=\sup\{|f_1(\varepsilon_1,\ldots,\varepsilon_{n^2})|,\ldots,|f_n(\varepsilon_1,\ldots,\varepsilon_{n^2})|:\sum|\varepsilon_{i,j}|\leq 1\}$$ 
is finite. 
Then from \eqref{detbl} we get
\begin{equation}\label{detbl2}
	|\det B_{\lambda}|\leq |\lambda|^n+M_0n|\lambda^{n}|,
\end{equation}
where we used triangle inequality and $|\lambda|^i\leq |\lambda|^n, i=0,1,\ldots,n$. Taking $M=(n+1)M_0$ we get the first inequality. To prove the second one we also use triangle inequality in \eqref{detbl} and obtain
\begin{equation}\label{detbl3}
	|\det B_{\lambda}|\geq |\lambda|^n-M_0n|\lambda|^{n-1}.
\end{equation}
Hence for $|\lambda|\geq 2M$ we have  
$$ |\det B_{\lambda}|\geq \frac12 |\lambda|^n,$$
which ends the proof of the statement.
$\hfil \square\ $

\begin{rmk}\label{rmk1}
With the same assumptions we can prove similarly the first inequality of Statement \ref{stm1} for the cofactors of matrix $B_{\lambda}$. Namely, if we dente cofactors of matrix $B_{\lambda}\in M_n(\mathbb{C})$ by $B_{i,j}$ then 
$|B_{i,j}|\leq M|\lambda|^{n-1}$.
\end{rmk}

\begin{stm}\label{stm2}
Let $A_{\lambda},B_{\lambda} \in M_n(\mathbb{C})$ such that $A_{\lambda}$ consist of Jordan blocks of egenvalue $\lambda$, where $|\lambda|$ is sufficiently large  and $\|B_{\lambda}-A_{\lambda}\|\leq 1$, where $\|A\|=\sum|a_{i,j}|$, then there exist constant $C$ such that  
$$\|B_{\lambda}\|\leq \frac{C}{|\lambda|}. $$ 
\end{stm}
{\it Proof} Without loss of generality we assume that $A_{\lambda}$ is a Jordan block. Using inversion formula to matrix $B_{\lambda}$ we get 
$$\|B_{\lambda}\|=|\det B_{\lambda}|^{-1}\sum|B_{i,j}|, $$
where $B_{i,j}$ are cofactors of matrix $B_{\lambda}$.
Using Statement \ref{stm1} and Remark \ref{rmk1} to estimate $|\det B_{\lambda}|$ and $|B_{i,j}|$ we obtain
$$\|B_{\lambda}\|\leq \frac{2}{|\lambda|^n}\cdot M|\lambda|^{n-1}=\frac{2M}{|\lambda|},$$
which ends the proof.
$\hfil \square\ $

\newpage
\renewcommand\refname{\Huge{References}}

\end{document}